\documentclass[final,leqno,onefignum,onetabnum]{siamltex1213}
\usepackage{amsmath}
\usepackage{epsfig}

\newcommand{\beq}{\begin{equation}}
\newcommand{\eeq}{\end{equation}}

\newcommand{\lf}{\left}
\newcommand{\rt}{\right}
\newcommand{\half}{\frac{1}{2}}

\newcommand{\tP}{\tilde{P}}

\newcommand{\oP}{\bar{P}}

\newcommand{\hP}{\hat{P}}
\newcommand{\hp}{\hat{p}}
\newcommand{\CP}{\mathcal{P}}

\newcommand{\al}{\alpha}
\newcommand{\de}{\delta}
\newcommand{\De}{\Delta}
\newcommand{\la}{\lambda}

\title{Joint Statistics of Random Walk on $Z^1$\\
and Accumulation 
of Visits} 

\author{J. K. Percus\footnotemark[1]\ \footnotemark[2]
\and O. E. Percus\footnotemark[1]}

\begin{document}

\renewcommand{\thefootnote}{\fnsymbol{footnote}}
\footnotetext[1]{Courant Institute of Mathematical Sciences, 
New York University (percus@cims.nyu.edu)}
\footnotetext[2]{Physics Department, New York University}

\maketitle
\slugger{mms}{xxxx}{xx}{x}{x--x}

\begin{abstract}
We obtain the joint distribution $P_N(X,K |Z)$ of the location $X$ of
a one-dimensional symmetric next neighbor random walk on the integer lattice,
and the number of times the walk has visited a specified site $Z$.  This
distribution has a simple form in terms of the one variable distribution
$p_{N'} (X')$, where $N'=N-K$ and $X'$ is a function of $X, \, K$, and $Z$.
The marginal distribution  of $X$ and $K$ are obtained, as well as their
diffusion scaling limits.
\end{abstract}



\pagestyle{myheadings}
\thispagestyle{plain}
\markboth{TEX PRODUCTION}{USING SIAM'S \LaTeX\ MACROS}

\section{Introduction}

Random walks in random environments are models for an enormous number of
biological, physical, ---  processes and when the walker alters the environment
during its sojourn, the resulting phenomenology is many-faceted.  This
paper addresses a forerunner, arguably, of such a situation, in which we are
content to obtain the distribution of the number of times $K$ the walker
has visited or is visiting a specified lattice site $Z$ after its starting,
but does not thereby alter the
properties of that site.  Extension to the reinforced random walk context
is our implicit direction, but will not yet be attempted.  However, we
will try to employ tools that do extend to such more structured problems.
Further, to keep our study as uncluttered as possible, we restrict our
attention to the blatantly simplest case of a symmetric random walk on 
the one-dimensional  integer lattice $Z^1$.

The specific question we will address is that of the correlation  between
the current location $X$ and the visit number $K$ at $Z$ since this will,
in work soon to be reported, inform the question of the evolution of the 
pattern of visits.  Our key result, Eq.~\eqref{3.17}, is that their joint 
distribution is expressible quite simply in terms of the elementary symmetric
random walk, thereby having a readily visualized asymptotic form, but that 
this form has a somewhat unusual structure.

\section{The symmetric Random Walk on $Z^1$}

Let us review briefly \cite{Ka, Fe} the underlying problem of the statistics of 
a symmetric random walk making next neighbor moves from its initial 
location at the origin, as represented by the probability distribution.
\beq\label{2.1}
p_0 (X) = \de_{X, 0}
\eeq
and ask for the $N^{th}$ step distribution $P_N(X)$, the probability 
that the walker has arrived at location $X$ after its $N^{th}$ random
move.  Clearly,
\beq\label{2.2}
\begin{split}
p_{N+1} (X) &= \half \lf(p_N(X-1)+ p_N (X+1)\rt) \\
&p_0(X) = \de_{X,0}, 
\end{split}
\eeq
or in terms of the generating function
\beq \label{2.3}
\hp \lf(\la, X \rt) =\sum^\infty_{N=0} \la^N \; p_N (x),
\eeq
we have instead
\beq\label{2.4}
\hp \lf(\la, X\rt) -p_0(X) = \frac{\la}{2}
\lf(\hp \,(\la, X+1) + \hp \,( \la, X-1)\rt).
\eeq
Eq.~\eqref{2.4}, with the convergence restriction $\lim_{|X|\to\infty}
\;\hp \lf(\la, X\rt)=0$, is solvable in standard fashion as
\beq\label{2.5}
\begin{split}
\hp \lf(\la, X\rt) &= \frac{\al^{-1}+\al}{\al^{-1} - \al}\;\al^{|X|} \\
&\text{ where }\; \al= \frac{1}{\la} - \sqrt{\frac{1}{\la^2}-1}
\\
&\qquad \lf(\text{from }\; \la= 2\al/(1+\al^2)\rt).
\end{split}
\eeq

On the other hand, the condition that the walker arrive at $X$ in $N$ steps
after starting at $0$ is equivalent to it having taken $\half \,
\lf(N+X\rt)$ steps of $+1$ and $\half\lf(N-X\rt)$ steps of $-1$, thereby 
identifying it with a binomial distribution of move probability $p=q=1/2$:
\beq\label{2.6}
p_N(X) = \frac{1}{2^N} \lf( 
\begin{array}{c}
N\\
\frac{1}{2} \lf(N + X\rt)
\end{array}
\rt),
\eeq
with the restriction of course that $X\equiv N (\!\!\mod 2)$.  The equivalence
of \eqref{2.5} and \eqref{2.6} will be useful in the sequel: $\hp$ is
useful for deconvolution  in step number (and has convenient algebraic
properties) whereas $p$ is directly  visualizable.  In fact, if we set
$Y=\half\, \lf(N+X\rt)$, $0\le Y \le N$, then  \eqref{2.6} takes
the standard form $2^{-N} \binom{N}{Y}$, allowing us to immediately quote
the $N$-asymptotic form $(x \equiv X/N^{1/2}$, but $\De \,X=2$)
\beq\label{2.7}
\CP(X) = \lim_{N\to\infty} \, N^{1/2}\,p_N\lf(x\sqrt{N}\rt)
=\frac{1}{\sqrt{2\pi}}\; e^{-x^{2}}.
\eeq

\section{Joint Distribution}

We now focus on a site $Z>0$ together with the number of times, $K$
it has been visited by the walker over the course of $N$ steps.  One can 
say that $K$ is reacting passively  to the passage of the walker, a 
situation that will be greatly expanded in future work, but we are at this
state simply interested in 
\begin{equation}\label{3.1}
P_N (X, K|Z),
\end{equation}
the probability that after $N$ steps from an initial
\begin{equation}\label{3.2}
P_0 (X,K|Z) = \de_{K,0}\, \de_{X,0}
\end{equation}
the walker is now at site $X$, and site $Z$ has been visited $K$ times.
The simplest way to deal with this situation is by constructing a 
generating function \cite{Wi} over $K$:
\beq\label{3.3}
\oP_N \lf(X,V|Z\rt) = \sum^\infty_{K=0} \,V^K\;P_N\lf(X,K |Z\rt),
\eeq
where $|V| \le 1$ allows convergence to be uniform in $N$.  The construction
\eqref{3.3} is clearly equivalent to inserting a multiplicative weight $V$
at each arrival at $Z$, and so satisfies the weighted version of 
Eq.~\eqref{2.2}:
\begin{align}
\oP_{N+1} \lf(X,V|Z\rt) &= \half \lf(1+(V-1) \, \de_{X,Z}\rt)
\lf(\oP_N \lf(X+1, V|Z\rt) + \oP_N \lf(X-1, V|Z\rt)\rt), \nonumber\\
&\qquad \oP_0 \lf(X,V|Z\rt) = \de_{X,0},\quad\text{ and } \label{3.4} \\
\hP\lf(\la, X, V|Z\rt) &\equiv \sum^\infty_{N=0} \la^N
\oP_N \lf(X,V|Z\rt).\nonumber
\end{align}
For further simplification, we again go over to the generating 
function over $N$, obtaining
\beq\label{3.5}
\begin{aligned}
\hP \lf(\la, X, V|Z\rt) &= \de_{X,0} + \frac{\la}{2} 
\lf( 1+(V-1)\,\de_{X,Z}\rt) \\
&\qquad\lf(\hP \lf(\la, X+1, V|Z\rt)+ \hP
\lf(\la, X-1, V|Z\rt) \rt),
\end{aligned}
\eeq
to be solved.

Eq.~\eqref{3.5} is a bit more complex than \eqref{2.4}, but yields to 
the same method of solution that was implicitly employed in obtaining
\eqref{2.5} from \eqref{2.4}.  We set up a covering of $Z^1$ by three
overlapping closed subspaces:
\beq\label{3.6}
\begin{aligned}
\text{[I]}\hspace{.85in} & X\in [-\infty, 0]\\
\text{[II]} \hspace{.85in} & X\in [0, Z]\\
\text{[III]} \hspace{.85in} & X\in [Z, \infty]
\end{aligned}
\eeq
and separate \eqref{3.5} into its actions on the corresponding open subspaces
and their boundaries
\begin{alignat}{2}
&X\to -\infty :  &\quad \hP &(\la, X,V|Z ) \to 0\nonumber \\
&X\in \text{(I)} :  &\hP &(\la, X,V|Z) \nonumber\\
&&  &\quad= \frac{\la}{2} \lf(\hP \lf(\la, X+1, V|Z\rt)
+ \hP \lf(\la, X-1, V|Z\rt) \rt)\nonumber\\
&X=0 : &\hP &(\la,0,V|Z) =1 + \nonumber\\
&&  &\quad \frac{\la}{2} \lf(\hP \lf(\la, 1, V|Z\rt) +
 \hP \lf(\la, -1, V|Z\rt) \rt)\nonumber\\
&X\in \text{(II)} :  &\hP &(\la, X,V|Z)  \label{3.7}\\
&&  &\quad = \frac{\la}{2} \lf(\hP \lf(\la, X+1, V|Z\rt) +
\hP\lf(\la, X-1, V|Z\rt)\rt)\nonumber\\
&X=Z :  &\hP &(\la, X,V|Z) \nonumber\\
&&  &\quad= \frac{\la V}{2} \lf(\hP \lf(\la, Z+1, V|Z\rt) +
\hP\lf(\la, Z-1, V|Z\rt)\rt)\nonumber\\
&X\in \text{(III)} :  &\hP &(\la, X,V|Z)  \nonumber\\
&&  &\quad= \frac{\la}{2} \lf(\hP \lf(\la, X+1, V|Z\rt) +
\hP\lf(\la, X-1, V|Z\rt)\rt)\nonumber \\
&X\to \infty :  &\hP &(\la, X,V|Z ) \to 0. \nonumber
\end{alignat}

Since the general solution of 
\beq\label{3.8}
a<X<b:\qquad f(X) = \frac{\la}{2} \lf(f(X+1) + f(X-1)\rt)
\eeq
is given by
\beq\label{3.9}
\begin{split}
a\le X \le b:\quad& f(X) =  A\, \al^{-X} + B\, \al^X\\
&\text{ where }\quad \al = \frac{1}{\la} - \lf(\frac{1}{\la^2} -1\rt)^{1/2},
\end{split}
\eeq
we have from \eqref{3.7}
\beq\label{3.10}
\begin{array}{llll}
X\in \text{[I]}&\!\!: &\quad &\hP^I \lf(\la, X, V|Z\rt) = A\, \al^{-X} \\[1mm]
X\in \text{[II]}&\!\!: &\quad &\hP^{II}\lf(\la, X, V|Z\rt) =
B\, \al^{X} + C\,\al^{-X} \\[1mm]
X\in \text{[III]}&\!\!: &\quad &\hP^{III} \lf(\la, X, V|Z\rt) = D\, \al^{X} 
\end{array}
\eeq
where $A, B, C, D$ are functions of $\la, V$, and $Z, \al$ of $\la$ alone.
The uniqueness of $\hP$ at $X=0$ and $X=Z$ then implies two equalities
from \eqref{3.10}, two from \eqref{3.7}:
\begin{alignat}{2}
&A &= &\;B+C\nonumber \\
&D\, \al^Z &=&\;B\,\al^Z + C\, \al^{-Z}\label{3.11}\\
&A &= &\;1 + \frac{\la}{2} \lf(B\,\al + C\,\al^{-1}+ A\,\al\rt) \nonumber\\
&D \, \al^Z &= &\;\frac{\la V}{2} \lf(B\,\al^{Z-1} + C \al^{1-Z}
+ D\al^{Z+1}\rt), \nonumber
\end{alignat}
immediately solved using $\la/2= \al/(1+\al^2)$ as 
\begin{alignat}{2}
&A&=& \frac{1+\al^2}{1-\al^2} \, \lf(1-\al^{2Z}\rt) + 
\frac{1+\al^2}{1+\al^2 -2 \,\al^2 \,V}\;V\, \al^{2Z}\nonumber\\
&B&=& \frac{1+\al^2}{1-\al^2}, \; C=- \frac{1+\al^2}{1-\al^2} \; \al^{2Z}
+ \frac{1+\al^2}{1+2^Z -2 \al^2 \,V} \, V\,\al^{2Z},\label{3.12}\\
&D&=& \frac{1+\al^2}{1+\al^2 - 2\al^2\,V} \, V. \nonumber
\end{alignat}
Substituting into \eqref{3.10} and using a unifying notation we conclude
that 
\beq\label{3.13}
\begin{split}
\hP \lf(\la, X, V|Z\rt) &= \frac{1+\al^2}{1-\al^2} 
\lf(\al^{|X|} - \al^Z \;\al^{|X-Z|}\rt)\\
&+ \frac{1+\al^2}{1+\al^2 -2\al^2\,V} \;V \,\al^Z \, \al^{|X-Z|}.
\end{split}
\eeq
Eq.~\eqref{3.13} can also be obtained by strictly combinatorial  reasoning,
but this approach does not extend to reinforced  walks.

Extracting $P_N(X, K|Z)$ from \eqref{3.13} means taking the coefficient of 
$\la^N\; V^K$, which is far from obvious.  But it can be done indirectly.
We first extract the coefficient of $V^K$:
\beq\label{3.14}
\begin{aligned}
&\tP \lf(\la, X, K|Z\rt) \equiv \text{ coef } V^K 
\text{ in  }\hP \lf(\la, X, V|Z\rt)\\
&\quad  = \frac{1+\al^2}{1-\al^2} \lf(\al^{|X|} - \al^Z \, \al^{|X-Z|}\rt)
\de_{K,0}  \\
&\hspace{1in}+ \lf(\frac{2 \al^2}{1+\al^2}\rt)^{K-1}
\, \al^Z\,\al^{|X-Z|} \lf(1-\de_{K,0}\rt),
\end{aligned}
\eeq
and then observe that
\beq\label{3.15}
\text{coef } (\la^N) \text{ in } \frac{1+\al^2}{1-\al^2}\; \al^{|Y|} = p_N(Y)
\eeq
whereas
\beq\label{3.16}
\begin{aligned}
&\text{coef } \lf(\la^N\rt) \lf(\frac{2 \al^2}{1+\al^2}\rt)^{K-1} 
\al^{|Y|} =  \text{coef } (\la^N) (\al \la)^{K-1} \, \al^{|Y|}\\
&=\text{coef } \lf(\la^{N+1-K}\rt) \al^{|Y| + K-1} = \text{ coef }
\lf(\la^{N+1-K}\rt) \frac{1-\al^2}{1+\al^2} \;
\frac{1+\al^2}{1-\al^2} \; \al^{|Y|+K-1} \\
&=\text{coef } \lf(\la^{N+1-K}\rt) \frac{\la}{2\al}
\;\frac{1+\al^2}{1-\al^2}
\lf(\al^{|Y|+K-1} - \al^{|Y|+K+1}\rt)\\
&=\text{coef } \lf(\la^{N-K}\rt) \half
\;\frac{1+\al^2}{1-\al^2}
\lf(\al^{|Y|+K-2} - \al^{|Y|+K}\rt)\\
&=\half \lf(p_{N-K} \lf( |Y|+K-2\rt)-p_{N-K} \lf(|Y|+K\rt)\rt)
\end{aligned}
\eeq
It follows from \eqref{3.14}, \eqref{3.15}, \eqref{3.16} that 
\beq\label{3.17}
\begin{split}
&P_N \lf(X, K|Z\rt) = \lf(p_N(X) - p_N \lf(Z+ |X-Z|\rt) \rt) \de_{K,0} \\
&\qquad+ \half \lf( p_{N-K} \lf(|X-Z| + Z+K-2\rt)
- p_{N-K} \lf(|X-Z| + Z + K\rt)\rt) \lf(1-\de_{K,0}\rt).
\end{split}
\eeq

\section{The $X$ and $K$ Marginals}

If we want just $P_N(X|Z)$, it is necessary to sum $P_N(X,K|Z)$ over $K$,
which is equivalent to setting $V=1$ in $\oP_N(X,V|Z)$.  A somewhat
simpler path is to first set $V=1$ in $\hP(\la, X, V|Z)$ of \eqref{3.13},
obtaining
\beq\label{4.1}
\hP \lf(\la, X, 1|Z\rt) = \frac{1+\al^2}{1-\al^2} \; \al^{|X|}\ ,
\eeq
and then convert from $\la$ to $N$:
\beq\label{4.2}
P_N \,(X|Z) = \oP_N \lf(X, 1|Z\rt) =  p_N(X)
= \lf(1/2^N\rt) \binom{N}{\half (N+X) }
\quad \text{if } \; N-X \equiv 0\pmod{2},
\eeq
the input symmetric random walk, as must be the case 

The $K$-marginal
\beq\label{4.3}
\CP_N \lf(K|Z\rt) = \sum^\infty_{X=-\infty} P_N \lf(X, K|Z\rt)
\eeq
is almost as simple, if less familiar.  Here, we only need to proceed 
directly: substitute \eqref{3.17} into \eqref{4.3} and carry out the 
summation over $X$.  In terms of the basic  cumulant
\beq\label{4.4}
F_N(X) = \sum^X_{Y=-\infty} p_N(Y)
\eeq
this is readily evaluated as
\beq\label{4.5}
\begin{split}
\CP_N (K|Z) &= F_N\,(Z-1) \, \de_{K,0}\\
&+\qquad \half \lf[ p_{N-K} \,(Z +K-2) + 2p_{N-K} \,(Z+K-1)
+ p_{N-K} \,(Z+K) \rt] \lf(1-\de_{K,0}\rt).
\end{split}
\eeq
or
\beq\label{4.6}
\CP_N (K|Z) = F_N(Z-1)\, \de_{K, 0} + \lf(1-\de_{K, 0}\rt)
\begin{cases}
1/2^{N-K} \binom{N-K}{\half \lf(N+Z-1\rt)} &\text{if }\; N-Z \equiv 1
\pmod{2}\\[2mm]
1/2^{N+1-K} 
\binom
{N+1-K}{\half \lf(N+Z\rt)}
&\text{if }\; N-Z \equiv 0
\pmod{2}
\end{cases}
\eeq

Eqn.~\eqref{3.17}, \eqref{4.2} and \eqref{4.5} solve the combinatorial problem
of determining the two-variable and two one-variable marginal distributions
of the random variables $X$ and $K$, but they don't present us with an 
immediate picture of what is going on .  That is traditionally done by finding
the $N$-dependence of various means and covariances, by looking for the 
distribution resulting for asymptotic $N$, or by both.  We adopt the third
strategy, which can be carried out in either order of its two components.  
To start, we may compute the low order moments.  Since these require only those
of the basic symmetric random walk, the results
\beq\label{4.7}
\begin{aligned}
E_N (X|Z) &= \sum^N_{X=-N} X \, P_N (X|Z) =0\\
E_N (X^2|Z) &= \sum^N_{X=-N} X^2 \, P_N (X|Z) \equiv N
\end{aligned}
\eeq
are immediate.  The corresponding $K$-moments are tedious to compute by 
hand (we have done so), but are eased by the use of Mathematica and we find
\beq\label{4.8}
\begin{aligned}
&E_N(K|Z) = \sum^N_{K=1} K\, P_N (K|Z)\\
& \qquad =
\begin{cases}
\sum^N_{K=1} \frac{K}{2^{N-K}} \binom{N-K}{\half \lf(N+Z-1\rt)}
\quad&\text{if} \quad N-Z\equiv 1 \pmod{2}\\[2mm]
\sum_{K=1}  \frac{K}{2^{N+1-K}} \binom{N+1-K}{\half \lf(N+Z\rt)} 
\quad&\text{if}\quad N-Z\equiv 0 \pmod{2}
\end{cases}\\
&\qquad 
\begin{array}{lll}
=&2^{1-N} \binom{N-1}{\half \lf(N+Z-1\rt)} 
&2 F1 \lf(2, \,\half \,(Z+1-N), \,1-N;2\rt)\\[2mm]
-\!\!&2(N+1) \binom{-1}{\half\lf(N+Z-1\rt)}
&2 F1 \lf(N+2,\, \half \,(Z+1+N),\,1+N,\,2\rt)\\
&& \text{if }\; N-Z \equiv 1\pmod{2}, \\[3mm]
=&2^{1-N} \binom{N}{\half \lf(N+Z\rt)} 
&2 F1 \lf(2, \,\half \,(Z-N), \,-N;Z\rt)\\[2mm] 
-\!\!&2(N+2) \binom{-1}{\half\lf(N+Z\rt)}
&2 F1 \lf(N+3,\, \half \,(Z+N+2),\,N+2;\,2\rt)\\  
&& \text{if }\; N-Z \equiv 10\pmod{2}, 
\end{array}
\end{aligned}
\eeq
in terms of Hypergeometric function $2F_1$.  Similarly, $E_N(K^2|Z)$
involves $3F2$.

\section{Asymptotics and Conclusion}

Although the explicit combinatorial results \eqref{3.17},
\eqref{4.7} do not contribute
to ready visualization, they do tell us about the nature of large $N$
asymptotics which, to put in a more positive light, can be regarded as 
conversion to a Brownian motion context. The simplification afforded by this
limit has as well the enormous  advantage of bringing out qualitative aspects
that may be concealed e.g., by analysis of various moments.  Let us see what 
this strategy reveals in our particular case.

In diffusion scaling, the transformation
\beq\label{5.1}
X=N^{1/2}\, x, \qquad \qquad K=N^{1/2}\, k
\eeq
is called for.  Since $\De k= 1/N^{1/2}$, but $X\equiv N \pmod{2}$
$\De x=2/N^{1/2}$, the probability in $(X, K)$ space is then transformed
as $N\to\infty$ to a probability density in $(x, k)$ space.  In addition,
the change of scale in $X$-space tells us that we should also scale $Z$ as
\beq\label{5.2}
Z=N^{1/2}\, z.
\eeq
The probability density in $(x,k)$ space then becomes
\beq\label{5.3}
\CP (x,k|z) =\lim_{N\to\infty} \:\frac{N}{2} \, P_N \lf(N^{1/2}\, x,
N^{1/2}\,k| N^{1/2}\,z\rt).
\eeq
We can now apply \eqref{2.7} to \eqref{3.17} (and use $\de_{K,0} \to N^{-1/2}
\; \de(k)$ with the convention that $\int^\infty_0 f(k) \, \de(k)\,dk=f(0)$) 
to obtain without difficulty
\beq\label{5.4}
\begin{split}
\CP (x,k|z) &= \frac{1}{\sqrt{2\pi}} \lf(e^{-\half\, x^2}-e^{-\half(z+|x-z|)^2}
\rt)\de(k) \\
&+ \half \frac{1}{\sqrt{2\pi}} \lf(k +z +|x-z|\rt)
e^{-\half \lf(k+z+|x-z|\rt)^2} \ ,
\end{split}
\eeq
divided of course into three regions, as was \eqref{3.6}.  The 
density maxima are shown (bold lines) in the accompanying figure.

\begin{figure}[h]
\centerline{\includegraphics[width=3in]{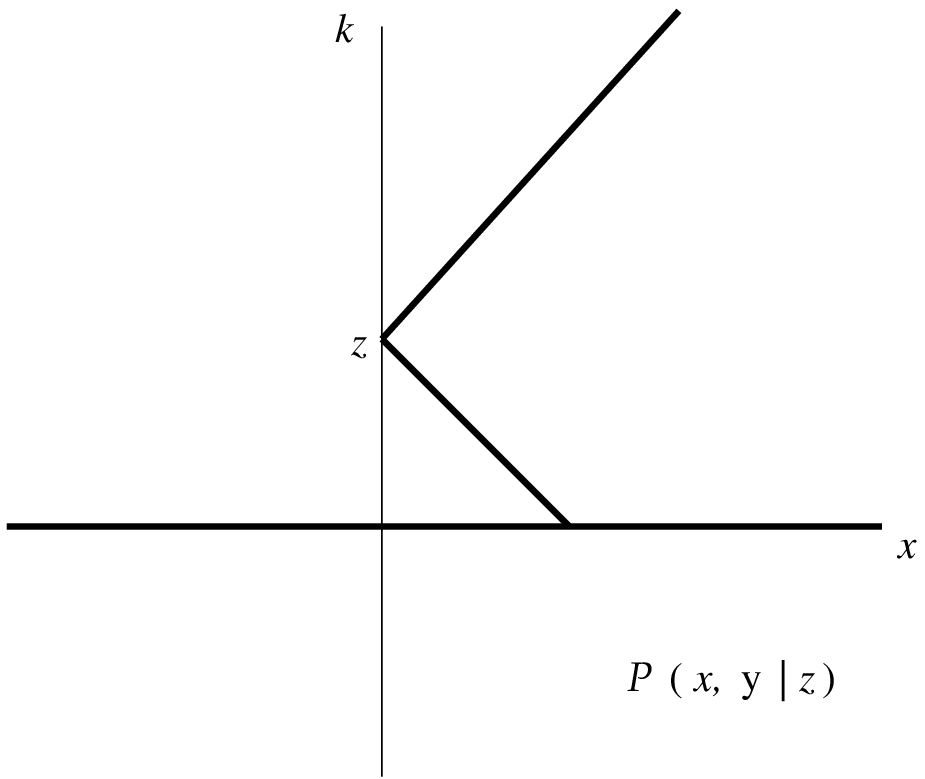}}
\end{figure}

Now the $x$-marginal is of course that of the basic symmetric Brownian
walk, but the $k$-marginal is available from either \eqref{4.6} or 
\eqref{5.4},
\beq\label{5.5}
P(k|z) = (2/ \pi)^{1/2} \, e^{-\half\, (k+z)^2} + C(z)\, \de(k)
\eeq
where $C(z)$ assures normalization.

It is clear that the first few joint moments, as in \eqref{4.7}, say very 
little about the structure of the joint distribution, that the asymptotic
form is both simpler and more informative.  (Note the $k+z$  dependence
of \eqref{5.5})
We will take advantage of this 
in work soon to be reported, in which the development of the full pattern of 
visits is in question.

\section*{Acknowledgment}

We are pleased to acknowledge the aid of an anonymous referee in improving
the final form of this paper.


\newpage

\begin{thebibliography}{00}
\bibitem{Ka}{\sc M. Kac}
{\em Random walk and the theory of Brownian motion}, 
Amer.\ Math Monthly \textbf{54} (1968), 369--391.
\bibitem{Fe}{\sc W. Feller}
{\em An Introduction to Probability Theory and Its Applications},
Vol.~1, J. Wiley, New York, 1968.
\bibitem{Wi}{\sc H. S. Wilf}
{\em Generating Functionology}, Academic Press, San Diego, 1990.
\end{thebibliography}
\end{document}